\documentclass{amsart}
\newtheorem{Thm}{Theorem}[section]

\newtheorem{Cor}[Thm]{Corollary}

\theoremstyle{remark}

\numberwithin{equation}{section}
\begin{document}

\title[First $L^p$-cohomology of groups with one end]
{The first $L^p$-cohomology of\\ some groups with one end}

\author[M. J. Puls]{Michael J. Puls}
\address{Department of Mathematics \\
Eastern Oregon University \\
One University Boulevard \\
LaGrande, OR 97850 \\
USA}
\email{mpuls@eou.edu}

\begin{abstract}
Let $p$ be a real number greater than one. In this paper we study the vanishing and nonvanishing of the first $L^p$-cohomology space of some groups that have one end. We also make a connection between the first $L^p$-cohomolgy space and the Floyd boundary of the Cayley graph of a group. We apply the result about Floyd boundaries to show that there exists a real number $p$ such that the first $L^p$-cohomology space of a nonelementary hyperbolic group does not vanish.
\end{abstract}


\date{November 3, 2006}
\maketitle

\section{Introduction}\label{Introduction}
In this paper $G$ will always be a finitely generated infinite group with identity 1 and symmetric generating set $S$. Let $\mathcal{F}(G)$ denote the set of all real-valued functions on $G$. Let $1 \leq p \in \mathbb{R}$ and set \[ D^p(G) = \{ f \in \mathcal{F}(G) \mid \sum_{g \in G} |f(gs^{-1}) - f(g) |^p < \infty \mbox{ for all }s \in S \}. \] The set $D^p(G)$ is known as the set of $p$-Dirichlet finite functions on $G$. Observe that the constant functions are in $D^p(G)$. We define a norm on $D^p(G)$ by \[ \parallel f \parallel_{D^p(G)} = \left( \left( \sum_{s \in S} \sum_{g \in G} |f(gs^{-1}) - f(g)|^p\right) + |f(1)|^p\right)^{\frac{1}{p}}. \] Under this norm $D^p(G)$ is a Banach space. We now define an equivalence relation on $D^p(G)$ by $f_1 \simeq f_2$ if and only if $f_1 - f_2$ is a constant function. Identify the constant functions by $\mathbb{R}$. Now $D^p(G) /\mathbb{R}$ is a Banach space under the norm induced from $D^p(G)$. That is, if $[f]$ is an equivalence class from $D^p(G) / \mathbb{R}$ then \[ \parallel [f]\parallel_{D^p(G)/\mathbb{R}} = \left( \sum_{s \in S} \sum_{g\in G} | f(gs^{-1}) - f(g) |^p \right)^{\frac{1}{p}}.\] We shall write $\parallel f \parallel_{D(p)}$ for $ \parallel [f] \parallel_{D^p(G)/\mathbb{R}}$. The norm for $D^p(G)$ and $D^p(G)/\mathbb{R}$ depends on the symmetric generating set $S$, but the underlying topology does not. If $A \subseteq D^p(G)$, then $(\overline{A})_{D^p(G)}$ will denote the closure of $A$ in $D^p(G)$. Similarly if $B \subseteq D^p(G)/\mathbb{R}$, then $(\overline{B})_{D(p)}$ will denote the closure of $B$ in $D^p(G)/\mathbb{R}$. Let $L^p(G)$ be the set that consists of functions on $G$ for which $\sum_{g\in G} |f(g)|^p$ is finite. Observe that $L^p(G)$ is contained in $D^p(G)/\mathbb{R}$. The main object of study in this paper is the space \[ \bar{H}_{(p)}^1 (G) = D^p(G)/ (\overline{L^p(G) \oplus \mathbb{R}})_{D(p)}. \] The space $\bar{H}_{(p)}^1 (G)$ is known as the first reduced $L^p$-cohomology space of $G$. This paper was inspired by the paper \cite{cartwoess}. 

It is well known that if $G$ has two ends then $\bar{H}_{(p)}^1 (G) = 0$ for $1 < p \in \mathbb{R}$. It is also well known that if $G$ has infinitely many ends then $\bar{H}_{(p)}^1 (G) \neq 0$ for $1 \leq p \in \mathbb{R}$, see \cite[Corollary 4.3]{Puls2} for a proof. A reasonable question to ask is: What can we say about $\bar{H}_{(p)}^1 (G)$ if $G$ has one end? It was shown in \cite[Corollary 3.6]{Puls2} that if $G$ has polynomial growth, then $\bar{H}_{(p)}^1 (G) = 0$ for $1 < p \in \mathbb{R}$. In \cite[Theorem 2]{BMV} it was shown that if $G$ is a properly discontinuous subgroup of isometries of a proper $CAT(-1)$ space with finite critical exponent and if the limit set of $G$ has at least three points, then $\bar{H}_{(p)}^1 (G) \neq 0$ for $p > \mbox{ max}\{ 1,\mbox{ critical exponent of } G \}$. Another result concerning groups with one end was given in \cite{Pansu} where it was shown that if $G$ is a co-compact lattice in $Sp(n,1)$, then $\bar{H}^1_{(p)} (G) \neq 0$ exactly for $p > 4n+2$.

Before we state our first result we need to define what it means for a Riemannian manifold to be rotationally symmetric. Let $M_n$ be a simply connected, $n$-dimensional Riemannian manifold with all sectional curvatures bounded above by a negative constant. Now fix a point on $M_n$ and use the exponential map at this point to transfer the polar coordinates on $\mathbb{R}^n$ to the manifold. So the Riemannian metric on $M_n$ can be written as $dx^2 = dr^2 + f(r)^2d\theta^2$ where $d\theta^2$ is the usual metric on the unit sphere $S^{n-1}, n \geq 2$. If the submanifolds $r = k$, where $k$ is a constant, are spheres of constant curvature then we shall say that $M_n$ is rotationally symmetric. In this paper we will prove:

\begin{Thm} \label{Riemannian}
Let $M_n$ be a complete, simply connected, $n$-dimensional Riemannian manifold with all sectional curvatures bounded above by a negative constant. Furthermore assume that $M_n$ is rotationally symmetric. Suppose that $G$ acts properly discontinuously on $M_n$ by isometries and that the action is cocompact and free. Then for $1 < p \leq n-1,\bar{H}_{(p)}^1 (G) = 0$.
\end{Thm}

What happens if $p > n-1$? Let $\mathcal{H}^n$ denote hyperbolic $n$-space. By combining Theorem 2 of \cite{BMV} with Theorem 1.6.1 of Nicholls \cite{Nicholls} we obtain the following: 
\begin{Thm} \label{hyperbolic}
Suppose that $G$ acts properly discontinuously on $\mathcal{H}^n$ by isometries and that the action is cocompact and free. If the limit set of $G$ has at least three points, then $\bar{H}_{(p)}^1 (G) \neq 0$ for $p > n-1$.
\end{Thm}

One of the hypothesis for Theorem \ref{hyperbolic} is that the limit set of $G$, which is a subset of the $(n-1)$-dimensional unit sphere, contain at least three points. Thus a possible first step in trying to determine whether $\bar{H}_{(p)}^1 (G)$ vanishes or does not vanish for groups with one end is to use a boundary for $G$ that is finer than the end boundary. One such boundary is the Floyd boundary. In Section \ref{Floydboundaries} we will prove
\begin{Thm} \label{Floyd}
Let $G$ be a finitely generated group and let $F$ be a Floyd admissible function on $G$. If the Floyd boundary of $G$ with respect to $F$ is nontrivial and if $\sum_{g\in G} (F(|g|))^p < \infty$, then $\bar{H}_{(p)}^1 (G) \neq 0$.
\end{Thm}
All concepts in Theorem \ref{Floyd} that are unfamiliar to the reader will be explained in Section \ref{Floydboundaries}. We will conclude Section \ref{Floydboundaries} by proving the following consequence, which appears to be known to Gromov (see pages 257-258 of \cite{Gromov}), of Theorem \ref{Floyd}.
\begin{Cor} \label{nonvanishhyper}
Let $G$ be a nonelementary hyperbolic group, then there exists a real number $p$ such that $\bar{H}^1_{(p)} (G) \neq 0.$
\end{Cor}

Let $f$ be an element of $\mathcal{F}(G)$ and let $g \in G$. Let $1< p \in \mathbb{R}$ and define
$$(\triangle_p f)(g) \colon = \sum_{s \in S} |f(gs^{-1}) - f(g)|^{p-2} ( f(gs^{-1}) - f(g)).$$
In the case $1 < p < 2$, we make the convention that $|f(gs^{-1}) - f(g)|^{p-2} (f(gs^{-1}) - f(g)) = 0$ if $f(gs^{-1}) = f(g)$. We shall say that $f$ is {\em $p$-harmonic} if $f \in D^p(G)$ and $\triangle_p f(g) = 0$ for all $g \in G$. Let $HD^p(G)$ be the set of $p$-harmonic functions on $G$. Observe that the constant functions are in $HD^p(G)$. If $p=2$, then $HD^2(G)$ is the linear space of harmonic functions on $G$ with finite energy. In general, $HD^p(G)$ is not a linear space if $p \neq 2$. A wealth of information about $p$-harmonic functions on graphs and manifolds can be found in \cite{H3, H2, H1}. The following decomposition theorem was proved in \cite[Theorem 3.5]{Puls2}.
\begin{Thm} \label{decomp}
Let $1 < p \in \mathbb{R}$ and suppose $(\overline{L^p(G)})_{D^p(G)} \neq D^p(G)$. Then for $f \in D^p(G)$, we can write $f = u + h$, where $u \in (\overline{L^p(G)})_{D^p(G)}$ and $h \in HD^p(G)$. This decomposition is unique up to a constant function.
\end{Thm}
It follows from the theorem that each nonzero class in $\bar{H}_{(p)}^1 (G)$ can be represented by a nonconstant function from $HD^p(G)$. This gives us the following:
\begin{Cor} \label{pharmonic}
Let $G$ be a finitely generated group
\begin{enumerate}
\item If $G$ satisfies the hypothesis of Theorem \ref{Riemannian}. Then $HD^p(G) = \mathbb{R}$ for $1 < p \leq n-1$.
\item If $G$ satisfies the hypothesis of Theorem \ref{hyperbolic}, then $HD^p(G)$ contains a nonconstant function for $p > n-1$.
\item If $G$ is a nonelementary hyperbolic group, then there exist a real number $p$ such that $HD^p(G)$ contains a nonconstant function.
\end{enumerate}
\end{Cor}

I would like to thank the referee for making many useful remarks that greatly improved the exposition of the paper. I would also like to thank Yaroslav Kopylov for some useful comments on a preliminary version of this paper. This work was supported by a grant from the research office at Eastern Oregon University. The author would like to thank the university for their kind support.

\section{Proof of Theorem 1.1} \label{Proofs}
In this section we will prove Theorem \ref{Riemannian}. We will begin by giving some definitions and other preliminaries needed for the proof of the theorem. Let $g \in G$ and let $f \in \mathcal{F}(G)$. Convolution of $f$ by $g-1$, denoted by $f \ast (g-1)$, is the function $(f \ast (g-1))(x) = f(xg^{-1}) - f(x)$ for $x \in G$. Observe that $f \ast (g-1) \in L^p(G)$ when $f \in D^p(G)$. The right translation of $f$ by $g$ is the function defined by $f_g(x) = f(xg^{-1})$. We will denote by $C_0(G)$ the set of those $f\in \mathcal{F}(G)$ for which the set $\{ g \mid |f(g)| > \epsilon \}$ is finite for each $\epsilon > 0$. The $L^p(G)$-norm for functions $f$ in $L^p(G)$ is denoted by $\parallel f \parallel_p$ and is given by $\parallel f \parallel_p^p = \sum_{g \in G} |f(g)|^p$.  

Let $M_n$ be a $n$-dimensional Riemannian manifold that satisfies the hypothesis of Theorem \ref{Riemannian}. The action of $G$ on $M_n$ will be denoted by $xg^{-1}$, where $x \in M_n$ and $g \in G$. The space $L^p(M_n)$ will consist of all real-valued functions on $M_n$ for which $\int_{M_n} |f(x)|^p dx < \infty$, where $1 \leq p \in \mathbb{R}$. We now proceed to prove the theorem. 

Suppose $\bar{H}_{(p)}^1 (G) \neq 0$ for some $p$ that satisfies $1 < p \leq n-1$. Then by the remark following Theorem \ref{decomp} there exists a nonconstant $p$-harmonic function $h$ in $HD^p(G)$ that represents a nonzero class in $\bar{H}_{(p)}^1(G)$. Define an affine isometric action of $G$ on $L^p(G)$ by $gf = f_g + h\ast (g-1)$. Let $(x_1, f_1)$ and $(x_2, f_2)$ be elements of the direct product $M_n \times L^p(G)$. We shall say that $(x_1, f_1)$ is related to $(x_2, f_2)$ if and only if there exists a $g \in G$ for which $x_2 = x_1g^{-1}$ and $f_2 = gf_1$. It is an easy exercise to show that this relation is an equivalence relation. Denote the quotient space of this equivalence relation by $M_n \times_G L^p(G)$. We now have a fibre bundle $M_n \times_G L^p(G) \stackrel{\pi}{\longrightarrow} M_n/G$, where $\pi$ denotes the projection map. Let $s$ be a smooth section of this bundle. Then $s(x) =(x,f) = (xg^{-1}, gf)$ where $\pi(x, f) = x$. We now define a smooth map $\hat{s}\colon M_n\rightarrow L^p(G)$ by $\hat{s}(x) = f$. Observe that  $\hat{s}(xg^{-1}) = g\hat{s} (x)$ since $(x, f) = (xg^{-1}, gf)$. We now define a real-valued function on $M_n$ by $f(x) \colon = (\hat{s}(x) + h)(1)$. If $g \in G$ and $x \in M_n$ then $f(xg^{-1}) = (g \hat{s}(x) + h)(1) = (\hat{s}(x)_g + h\ast (g-1) + h)(1) = \hat{s}(x) (g^{-1}) + h(g^{-1}) = (\hat{s}(x) + h)(g^{-1})$. Since $df(x) = d\hat{s}(x)(1)$ it now follows that $df(xg^{-1}) = d\hat{s}(x)(g^{-1})$, where $df$ is the differential of $f$.  Due to the compactness of $M_n/G$ there exists a constant $C$ such that $\sum_{g\in G} | d\hat{s}(x)(g^{-1})|^p = \parallel d\hat{s}(x) \parallel_p^p < C$ for all $x \in M_n$. Thus $\int_{M_n} |df(x)|^p dx = \int_{M_n/G} \sum_{g\in G} |df(xg^{-1})|^p dV \leq C(\mbox{ volume}(M_n/G)).$ Hence $df \in L^p(M_n)$ by the compactness of $M_n/G$. Using the canonical identification of $df$ with $\nabla f$, page 160 of \cite{Darling}, where $\nabla f$ is the gradient of $f$, we see that $\nabla f \in L^p(M_n)$. By \cite[Theorem 5.8]{Strichartz} there exists a constant $c$ such that $f-c \in L^p(M_n)$. Thus $\int_{M_n/G} \sum_{g \in G} |(f-c)(xg^{-1})|^p dV = \int_{M_n} |(f-c)(x)|^p dx < \infty$. Hence, $\sum_{g \in G} |(f-c)(xg^{-1})|^p < \infty$ for a fixed $x \in M_n$. Consequently $\hat{s}(x) + h - c \in L^p(G)$. Thus the $p$-harmonic function $h-c \in L^p(G) \subseteq C_0(G)$ since $\hat{s}(x) \in L^p(G)$. Lemma 6.1 of \cite{Puls2} tells us that $h-c = 0$ on $G$, contradicting the fact that $h$ is nonconstant. Therefore, $\bar{H}_{(p)}^1(G) = 0$ for $1 < p \leq n-1$. This concludes the proof of Theorem \ref{Riemannian}

\section{Floyd Boundaries} \label{Floydboundaries}
Let $(X, S)$ be the Cayley graph of $G$ with respect to the generating set $S$. Thus the vertices of $(X,S)$ are the elements of $G$, and $g_1, g_2 \in G$ are joined by an edge if and only if $g_1 = g_2s^{\pm 1}$ for some generator $s$. For the rest of this paper we will denote $(X,S)$ by $X$. We can make $X$ into a metric space by assigning length one to each edge, and defining the distance $d_s(g, h)$ between any two vertices $g,h$ in $X$ to be the length of the shortest path between $g$ and $h$. The metric $d_s$ on $X$ is known as the word metric. For the rest of this paper we will drop the use of the subscript $s$ and $d(x,y)$ will always denote the distance between $x$ and $y$ in the word metric. We will denote $d(1,g)$ by $|g|$ for $g\in G$. If $A$ is a set of vertices from $X$, then $|A| = \inf_{a \in A} d(1,a)$. Let $F$ be a function from the natural numbers $\mathbb{N}$ into the positive real numbers $\mathbb{R}^+$. We shall say that $F$ is a Floyd admissible function if it is monotonically decreasing, summable and for which there is a positive constant $L$ that satisfies $F(n+1) \leq F(n) \leq L F(n+1)$ for all $n \in \mathbb{N}$. We will now show how to construct a Floyd boundary for $X$ with respect to $F$. First we use $F$ to define a new metric on $X$. The new length of an edge joining $g$ and $h$ is $F(|\{ g, h\}|)$. Let $\alpha = \{ g_i\}$ be a path in $X$. The length $L_F$ of $\alpha$ is given by $\sum_{i} F(|\{ g_i, g_{i+1}\} | )$ and the new distance between $x$ and $y$ in $X$ is $d_F (x, y) := \inf_{\alpha} L_F(\alpha)$, where the infimum is taken over all paths $\alpha$ connecting $x$ and $y$. It is straight forward to verify that $d_F$ is a metric on $X$. Let $(\overline{X}^F, \bar{d}_F)$ denote the completion of $(X, d_F)$ in the sense of metric spaces. The {\em Floyd boundary} of $X$ is the set $\partial_F X = \overline{X}^F \backslash X$. We shall say that $\partial_F X$ is trivial if it consists of only $0, 1$ or $2$ points. Lots of information about Floyd boundaries can be found in \cite{Karlsson2, Karlsson3, Karlsson1}. If $A$ is a set, then the cardinality of $A$ will be denoted by $\#(A)$.

We now prove Theorem \ref{Floyd}. Let $f$ be a continuous functions from $\overline{X}^F$ into $\mathbb{R}$ that satisfies a Lipshitz condition. Thus

\begin{equation*}
\begin{split}
\sum_{g \in G} \sum_{s \in S} |f(gs^{-1}) - f(g)|^p & \leq \sum_{g\in G} \sum_{s\in S} C(d_F(gs^{-1}, g))^p \\
  & \leq \sum_{g\in G} \sum_{ s \in S} C (F(|g|))^p \\
   & = \#(S)C \sum_{g\in G} (F(|g|))^p < \infty.
\end{split}
\end{equation*}
So $f$ restricted to $G$ is an element of $D^p(G)$. Let $\xi \in \partial X_F$. Define a real-valued function on $\overline{X}^F$ by $f(x) \colon = \bar{d}_F (x, \xi).$ Now $f \in D^p(G)$ since it satisfies a Lipshitz's condition.  Since $\bar{d}_F$ is a metric $f$ is nonconstant on $\partial X_F$. By continuity of $f$ it follows that $f(|g|)$ does not tend towards a constant number as $|g|$ goes to infinity in $X$. Thus $f$ is not an element of $L^p(G) \oplus \mathbb{R}$. Also $G$ is nonamenable since $\partial X_F$ contains more than two points, \cite[Corollary 2]{Karlsson3}. Hence $L^p(G) \oplus \mathbb{R}$ is closed in $D^p(G)$ \cite[Corollary 1]{guichardet}, also see\cite[Theorem 4.1]{Puls2} for a proof. Therefore $f$ represents a nonzero class in $\overline{H}^1_{(p)} (G)$ and the proof of Theorem \ref{Floyd} is now complete.

We will now apply Theorem \ref{Floyd} to a class of hyperbolic groups. For the rest of this section assume that $G$ is a hyperbolic group with hyperbolic constant $\delta$. Let $X$ be the Cayley graph of $G$. The Gromov inner product with basepoint 1 in $X$ is defined to be
$$ (x\mid y ) = \frac{1}{2} ( |x| + |y| - d(x,y))$$
where $x$ and $y$ are vertices in $X$. Let $(x_n)$ be a sequence in $X$. We shall say that $(x_n)$ converges to infinity if $\lim_{n,m \rightarrow \infty} (x_n \mid x_m) = \infty$. Let $S_{\infty} (X)$ be the set of all sequences on $X$ which converge to infinity. We shall also say that two sequences, $(x_n)$ and $(y_n)$ in $S_{\infty} (X)$ are related if and only if $\lim_{n \rightarrow \infty} (x_n \mid y_n) = \infty$. This relation is an equivalence relation since $G$ is hyperbolic. The sequential boundary of $G$, denoted by $\partial G$, is the set of equivalence classes of sequences under the above relation. A hyperbolic group is nonelementary if there are more than two elements in $\partial G$.

Choose $a > 0$ such that $e^{3\delta a} - 1 < \sqrt{2} - 1$. Define a Floyd admissible function $F$ from $\mathbb{N}$ into $\mathbb{R}^+$ by $F(n) = e^{-an}$. Let $g\in G$ and $s \in S$. Then $g$ and $gs^{-1}$ are neighbors in $X$ and $F( | \{ g, gs^{-1}\}|) = e^{-a ( g \mid gs^{-1})}$. Thus for $x$ and $y$ in $X$, $d_F (x,y) = \inf \{ \sum_{i=1}^n e^{-a(x_i \mid x_{i+1})} \mid n \geq 1, x=x_0, x_1, x_2, \dots, x_n = y \in X\}$. Where, $x_0, x_1, \dots , x_n$ is a path from $x$ to $y$ in $X$. Now let $(x_n)$ and $(y_n)$ be sequences in $X$. By Proposition 22.8 of \cite{woess} we have the following inequality
$$ (3 - 2e^{3\delta a})e^{-a(x_n \mid y_n)} \leq d_F(x_n, y_n) \leq e^{-a (x_n \mid y_n )}.$$
Thus $\lim_{n \rightarrow \infty} d_F (x_n, y_n) = 0$ if and only if $\lim_{n \rightarrow \infty} ( x_n \mid y_n) = \infty$. Hence, the cardinality of $\partial X_F$ equals the cardinality of $\partial G$.

Since $G$ is finitely generated it has at most exponential growth. Consequently, there exists a real number $p$ such that
\begin{equation*}
\begin{split}
\sum_{g \in G} |F(|g|)|^p  & = \sum_{g \in G} e^{-a|g|p} \\
                            & \leq 2k \sum_{n = 1}^{\infty} (2k-1)^{n-1} e^{-anp} \\
                           & < \infty,
\end{split}
\end{equation*}
where $2k$ is the cardinality of $S$. Now apply Theorem \ref{Floyd} to obtain Corollary \ref{nonvanishhyper}.

\bibliographystyle{plain}
\bibliography{revoneendlp}

\begin{thebibliography}{10}

\bibitem{BMV}
Marc Bourdon, Florian Martin, and Alain Valette.
\newblock Vanishing and non-vanishing for the first ${L}^p$-cohomology of
  groups.
\newblock {\em Comment. Math. Helv.}, 80:377--389, 2005.

\bibitem{cartwoess}
Donald~I. Cartwright and Wolfgang Woess.
\newblock Infinite graphs with nonconstant {D}irichlet finite harmonic
  functions.
\newblock {\em SIAM J. Discrete Math.}, 5(3):380--385, 1992.

\bibitem{Darling}
R.~W.~R. Darling.
\newblock {\em Differential forms and connections}.
\newblock Cambridge University Press, Cambridge, 1994.

\bibitem{Gromov}
M.~Gromov.
\newblock Asymptotic invariants of infinite groups.
\newblock In {\em Geometric group theory, Vol.\ 2 (Sussex, 1991)}, volume 182
  of {\em London Math. Soc. Lecture Note Ser.}, pages 1--295. Cambridge Univ.
  Press, Cambridge, 1993.

\bibitem{guichardet}
A.~Guichardet.
\newblock \'{E}tude de la {$l$}-cohomologie et de la topologie du dual pour les
  groupes de {L}ie \`a radical ab\'elien.
\newblock {\em Math. Ann.}, 228(3):215--232, 1977.

\bibitem{H3}
Ilkka Holopainen.
\newblock Rough isometries and {$p$}-harmonic functions with finite {D}irichlet
  integral.
\newblock {\em Rev. Mat. Iberoamericana}, 10(1):143--176, 1994.

\bibitem{H2}
Ilkka Holopainen and Paolo~M. Soardi.
\newblock {$p$}-harmonic functions on graphs and manifolds.
\newblock {\em Manuscripta Math.}, 94(1):95--110, 1997.

\bibitem{H1}
Ilkka Holopainen and Paolo~M. Soardi.
\newblock A strong {L}iouville theorem for $p$-harmonic functions on graphs.
\newblock {\em Ann. Acad. Sci. Fenn. Math.}, 22(1):205--226, 1997.

\bibitem{Karlsson2}
Anders Karlsson.
\newblock Boundaries and random walks on finitely generated infinite groups.
\newblock {\em Ark. Mat.}, 41(2):295--306, 2003.

\bibitem{Karlsson3}
Anders Karlsson.
\newblock Free subgroups of groups with nontrivial {F}loyd boundary.
\newblock {\em Comm. Algebra}, 31(11):5361--5376, 2003.

\bibitem{Karlsson1}
Anders Karlsson.
\newblock Some remarks concerning harmonic functions on homogeneous graphs.
\newblock In {\em Discrete random walks (Paris, 2003)}, Discrete Math. Theor.
  Comput. Sci. Proc., AC, pages 137--144 (electronic). Assoc. Discrete Math.
  Theor. Comput. Sci., Nancy, 2003.

\bibitem{Nicholls}
Peter~J. Nicholls.
\newblock {\em The ergodic theory of discrete groups}, volume 143 of {\em
  London Mathematical Society Lecture Note Series}.
\newblock Cambridge University Press, Cambridge, 1989.

\bibitem{Pansu}
Pierre Pansu.
\newblock Cohomologie {$L\sp p$} des vari\'et\'es \`a courbure n\'egative, cas
  du degr\'e {$1$}.
\newblock {\em Rend. Sem. Mat. Univ. Politec. Torino}, (Special Issue):95--120
  (1990), 1989.
\newblock Conference on Partial Differential Equations and Geometry (Torino,
  1988).

\bibitem{Puls2}
Michael~J. Puls.
\newblock The first ${L}^p$-cohomology of some finitely generated groups and
  $p$-harmonic functions.
\newblock {\em J. Funct. Anal.}, 237(2):391--401, 2006.

\bibitem{Strichartz}
Robert~S. Strichartz.
\newblock Analysis of the {L}aplacian on the complete {R}iemannian manifold.
\newblock {\em J. Funct. Anal.}, 52(1):48--79, 1983.

\bibitem{woess}
Wolfgang Woess.
\newblock {\em Random walks on infinite graphs and groups}, volume 138 of {\em
  Cambridge Tracts in Mathematics}.
\newblock Cambridge University Press, Cambridge, 2000.

\end{thebibliography}
\end{document}